\documentclass{article}

\title{\LARGE \textbf{On Longest Cycle $C$ of a Graph $G$ \\via Structures of $G-C$}}
\author{Zh.G. Nikoghosyan\footnote{G.G. Nicoghossian (up to 1997)}\\
Institute for Informatics and Automation Problems\\ National Academy of Sciences\\
P. Sevak 1, Yerevan 0014, Armenia\\ E-mail: zhora@ipia.sci.am}

\begin{document}

\maketitle

\begin{abstract}
Two sharp lower bounds for the length of a longest cycle $C$ of a graph $G$ are presented in terms of the lengths of a longest path and a longest cycle of $G-C$, denoted by $\overline{p}$ and $\overline{c}$, respectively, combined with minimum degree $\delta$: (1) $|C|\geq(\overline{p}+2)(\delta-\overline{p})$ and (2) $|C|\geq(\overline{c}+1)(\delta-\overline{c}+1)$.\\


\end{abstract}

\section{Introduction}

We consider only finite undirected graphs without loops or multiple edges. A good reference for any undefined terms is [1]. The set of vertices of a graph $G$ is denoted by $V(G)$ or just $V$; the set of edges by $E(G)$ or just $E$. For $S$ a subset of $V(G)$, we denote by $G-S$ the maximum subgraph of $G$ with vertex set $V(G)-S$. For a subgraph $H$ of $G$ we use $G-H$ short for $G-V(H)$. 

Paths and cycles in a graph $G$ are considered as subgraphs of $G$. If $Q$ is a path or a cycle then the length of $Q$, denoted by $|Q|$, is $|E(Q)|$. For $Q$ a path, we denote $|Q|=-1$ if and only if $V(Q)=\emptyset$. Throughout the paper each vertex and edge can be interpreted as cycles of lengths 1 and 2, respectively. The length of a longest cycle of $G$ is called a circumference. 

Almost all lower bounds for the circumference are based on a standard procedure: choose any initial cycle $C_0$  in a graph $G$ and try to enlarge it via structures of $G-C_0$ and connections between $C_0$ and $G-C_0$. This can be realized by deleting some segment of $C_0$ of the type $x\overrightarrow{C}_0y$ and adding appropriate $(x,y)$-paths passing through $G-C_0$. In practice, mainly the maximal paths of $G-C_0$ are used (combined with connectivity conditions) as an optimal structure in $G-C_0$ to enlarge $C_0$.   

In view of these motivations, for $C$ a longest cycle of $G$, the length of a longest path of $G-C$ or maybe another convenient parameter of $G-C$, would be the frequently appeared parameter incorporated into various lower bounds for the circumference. However, in practice we do not meet these expected parameters appeared so frequently. Instead, the degree and connectivity conditions are used in the majority of results. 

In this paper we present the first two lower bounds for the length of a longest cycle $C$ based on two parameters of $G-C$, namely the lengths of a longest path and a longest cycle of $G-C$, denoted by $\overline{p}$ and $\overline{c}$, respectively, combined with minimum degree $\delta$.\\

\noindent\textbf{Theorem 1}. For $C$ a longest cycle of a graph, $|C|\geq(\overline{p}+2)(\delta-\overline{p})$. \\

\noindent\textbf{Theorem 2}. For $C$ a longest cycle of a graph, $|C|\geq(\overline{c}+1)(\delta-\overline{c}+1)$.\\

The limit example $(\kappa+1)K_{\delta-\kappa+1}+K_{\kappa}$ shows that Theorems 1 and 2 are sharp. The first preprint versions of Theorems 1 and 2 can be found in [2] and [3], appeared still in 1998 and 2000, respectively. This preprint version aims to combine these two results in a united terminology and format. 

\section{Terminology}

An $(x,y)$-path is a path with endvertices $x$ and $y$. Given an $(x,y)$-path $L$ of $G$, we denote by $\overrightarrow{L}$ the path $L$ with an orientation from $x$ to $y$. If $u,v\in V(L)$, then $u\overrightarrow{L}v$ denotes the consecutive vertices on $\overrightarrow{L}$ from $u$ to $v$ in the direction specified by $\overrightarrow{L}$. The same vertices, in reverse order, are given by $v\overleftarrow{L}u$. For $\overrightarrow{L}=x\overrightarrow{L}y$ and $u\in V(L)$, let $u^+(\overrightarrow{L})$ (or just $u^+$) denotes the successor of $u$ $(u\neq y)$ on $\overrightarrow{L}$, and $u^-$ denotes its predecessor $(u\neq x)$. If $A\subseteq V(L)-y$ and $B\subseteq V(L)-x$, then we denote $A^+=\{v^+|v\in A\}$ and $B^-=\{v^-|v\in B\}$. Similar notation is used for cycles. If $Q$ is a cycle and $u\in V(Q)$, then $u\overrightarrow{Q}u=u$. For $v\in V$, put $N(v)=\{u\in V|uv\in E\}$ and $d(v)=|N(v)|$.

\section{Special definitions}

For the remainder of this section, let a subgraph $H$ of a graph $G$ and a path (or a cycle) $\overrightarrow{M}$ in $G-H$ be fixed and let $u_1,...,u_m$ be the vertices of $M$ occuring on $\overrightarrow{M}$ in a consecutive order. \\

\noindent\textbf{Definition 1} $\{M_H$-spreading; $\overrightarrow{\Upsilon}(u)$; $\dot{u}$; $\ddot{u}\}$. An $M_H$-spreading $\Upsilon$ is a family of pairwise disjoint paths $\overrightarrow{\Upsilon}(u_1),...,\overrightarrow{\Upsilon}(u_m)$ in $G-H$ with $\overrightarrow{\Upsilon}(u_i)=u_i\overrightarrow{\Upsilon}(u_i)\ddot{u}_i$ $(i=1,...,m)$. If $u\neq \ddot{u}$ for some $\overrightarrow{\Upsilon}(u)$, then we use $\dot{u}$ to denote the successor of $u$ along $\overrightarrow{\Upsilon}(u)$.\\  

\noindent\textbf{Definition 2} $\{\Phi_u; \varphi_u; \Psi_u; \psi_u\}$. Let $\Upsilon$ be any $M_H$-spreading. For each $u\in V(M)$, put 
$$
\Phi_u=N(\ddot{u})\cap V(\Upsilon), \quad \varphi_u=|\Phi_u|,\qquad\qquad\qquad\qquad
$$
$$
\Psi_u=N(\ddot{u})\cap V(H), \quad \psi_u=|\Psi_u|.\qquad\qquad\qquad\qquad
$$
\\
\noindent\textbf{Definition 3} $\{U_0; \overline{U}_0; U_1; U^*\}$. For $\Upsilon$ an $M_H$-spreading, put
$$
U_0=\{u\in V(M)|u=\ddot{u}\}, \quad  \overline{U}_0=V(M)-U_0,\qquad\qquad\qquad\quad\enspace
$$
$$
U^*=\{u\in \overline{U}_0|\Phi_u \subseteq V(\Upsilon(u))\}, \quad U_1=V(M)-(U_0\cup U^*).\qquad
$$
\\
\noindent\textbf{Definition 4} $\{(U_0)$-minimal and $(U_0,U^*)$-minimal $M_H$-spreadings$\}$. An $M_H$-spreading $\Upsilon$ is said to be $(U_0)$-minimal, if it is chosen such that $|U_0|$ is minimum. A $(U_0)$-minimal $M_H$-spreading is said to be $(U_0,U^*)$-minimal if it is chosen such that $|U^*|$ is minimum.\\

\noindent\textbf{Definition 5} $\{B_u; B^*_u; b_u; b^*_u\}$. For $\Upsilon$ an $M_H$-spreading and $u\in V(M)$, set $B_u=\{v\in U_0|v\dot{u}\in E\}$ and $b_u=|B_u|$. Further, for each $u\in U_0$, set $B^*_u=\{v\in \overline{U}_0|u\dot{v}\in E\}$ and $b^*_u=|B^*_u|$.

\section{Preliminaries}

\textbf{Lemma 1}. Let $C$ be a longest cycle in a graph $G$ and $M$ a path in $G-C$. Let $\overrightarrow{L}_1,...,\overrightarrow{L}_r$ be vertex disjoint paths in $G-C$ with $\overrightarrow{L}_i=v_{i}\overrightarrow{L}_iw_i$ $(i=1,...,r)$ having only $v_1,...,v_r$ in common with $M$. Then
$$
|C|\ge\sum_{i=1}^r|Z_i|+\Big|\bigcup_{i=1}^rZ_i\Big|,\qquad\qquad\qquad
$$
where $Z_i=N(w_i)\cap V(C)$ $(i=1,...,r)$.  \\

\noindent\textbf{Lemma 2}. Let $H$ be any subgraph of a graph $G$ and $M$ a longest cycle in $G-H$ with a $(U_0)$-minimal $M_H$-spreading $\Upsilon$. Then for each $u\in U_1$, $|M|\geq\varphi_u+b_u+1$. \\

\noindent\textbf{Lemma 3}. Let $H$ be any subgraph of a graph $G$ and $L$ a longest path in $G-H$ with a $(U_0)$-minimal $L_H$-spreading $\Upsilon$. Then for each $u\in \overline{U}_0$, $|L|\geq\varphi_u+b_u$.

\section{Proofs}

\textbf{Proof of lemma 1}. Assume without loss of generality (w.l.o.g.) that $v_i=w_i$ $(i=1,...,r)$ since otherwise, we can use the same arguments. If $\bigcup_{i=1}^r Z_i=\emptyset$, then there is nothing to prove. Let $\bigcup_{i=1}^r Z_i\neq\emptyset$ and let  $\xi_1,...,\xi_t$  be the elements of $\bigcup_{i=1}^r Z_i$  occurring on $\overrightarrow{C}$ in a consecutive order. Set 
\[
 F_i=N(\xi_i)\cap\lbrace w_1,...,w_r\rbrace\quad (i=1,...,t).\qquad\qquad\qquad\qquad\quad\enspace
\]

Suppose first that $t=1$. If $|F_1|=1$, then $\sum_{i=1}^r |Z_i|=|\bigcup_{i=1}^r Z_i|=1$ and the result follows from $|C|\geq 2$ immediately. If $|F_1|\geq 2$ then choosing a largest segment $u\overrightarrow{M}v$ on $M$ with $u,v\in F_1$, we get a cycle $C^{\prime}=\xi_1 u\overrightarrow{M}v\xi_1$ satisfying
\[
 |C|\geq|C^{\prime}|\geq\sum_{i=1}^r |Z_i|+1=\sum_{i=1}^r |Z_i|+\Big|\bigcup_{i=1}^r Z_i\Big|.\qquad\qquad\qquad\enspace
\]

Now assume $t\geq 2$. Putting  $f(\xi_i)=|\xi_i\overrightarrow{C}\xi_{i+1}|$ (indices mod $t$) for each $i\in\{1,...,t\}$, it is easy to see that 
$$
|C|=\sum_{i=1}^t f(\xi_i),\quad \sum_{i=1}^t |F_i|=\sum_{i=1}^r |Z_i|,\quad t=\Big|\bigcup_{i=1}^r Z_i\Big|.\qquad\qquad\eqno(1)
$$

For each $i\in \lbrace 1,...,t\rbrace$, let $x_i \overrightarrow{M}y_i$ be the largest segment on $\overrightarrow{M}$ with $x_i,y_i\in F_i\cup F_{i+1}$ (indices mod $t$). Now we need to show that $f(\xi_i)\ge (|F_i|+|F_{i+1}|+2)/2$. Indeed, if $x_i\in F_i$ and $y_i\in F_{i+1}$ then $f(\xi_i)\geq|\xi_i x_i \overrightarrow{M}y_i \xi_{i+1}|$ (since $C$ is extreme) implying that 
\[
f(\xi_i)\geq \max\lbrace|F_i|,|F_{i+1}|\rbrace +1\geq\frac{1}{2}(|F_i|+|F_{i+1}|+2).\qquad\qquad
\]

The same inequality holds from $f(\xi_i)\geq |\xi_i y_i \overleftarrow{M}x_i \xi_{i+1}|$ if $x_i\in F_{i+1}$ and $y_i\in F_i$, by a similar argument. Now suppose that either $x_i,y_i\in F_i$ or $x_i,y_i\in F_{i+1}$. Assume w.l.o.g. that $x_i,y_i\in F_i$. In addition, we have $x_i,y_i\notin F_{i+1}$, since otherwise we are in the previous case. Let $x_{i}^{\prime} \overrightarrow{M}y_{i}^{\prime}$ be the largest segment on $\overrightarrow{M}$ with $x_i ^{\prime},y_i ^{\prime} \in F_{i+1}$. If $|x_i \overrightarrow{M}x_i ^\prime |\ge(|F_i|-|F_{i+1}|)/2$ then $f(\xi_i)\ge|\xi_i x_i \overrightarrow{M}y_i^\prime \xi_{i+1}|$ and hence
$$
f(\xi_i)\ge\frac{1}{2}(|F_i|-|F_{i+1}|)+|F_{i+1}|+1=\frac{1}{2}(|F_i|+|F_{i+1}|+2).\enspace
$$

Finally, if $|x_i \overrightarrow{M}x_i^\prime |\le(|F_i|-|F_{i+1}|-1)/2$, then 
$$
f(\xi_i)\geq|\xi_iy_i\overleftarrow{M}x^{\prime}_i\xi_{i+1}|=|x_i^\prime\overrightarrow{M}y_i|+2=|x_i\overrightarrow{M}y_i|-|x_i\overrightarrow{M}x_i^\prime|+2
$$
$$
\qquad \qquad \ge|F_i|-1-\frac{1}{2}(|F_i|-|F_{i+1}|-1)+2>\frac{1}{2}(|F_i|+|F_{i+1}|+2).
$$

So, $f(\xi_i)\ge (|F_i|+|F_{i+1}|+2)/2 \quad (i=1,...,t)$ in any case. Therefore,
$$
\sum_{i=1} ^t f(\xi_i)\ge \sum_{i=1} ^t \frac{1}{2}(|F_i|+|F_{i+1}|+2)=\sum_{i=1}^t |F_i|+t\qquad\qquad\quad
$$
and the result follows from (1).   \quad  $ \Delta$\\

\noindent \textbf{Proof of lemma 2}. Let $u\in U_1$. For each $x\in V(M)$, put $A_u(x)=(\Phi_u\cup B_u)\cap V(\Upsilon(x))$. By the definition, 
$$
|\Phi_u\cup B_u|=\sum_{x\in V(M)}|A_u(x)|. \eqno(2)
$$

If $A_u(x)\neq\emptyset$ for some $x\in V(M)$, then we choose a vertex $\rho_u(x)$ in $A_u(x)$ such that $|x\overrightarrow{\Upsilon}(x)\rho_u(x)|$ is maximum. By the definition, $\rho_u(u)=(\ddot{u})^-$. Put $\overline{\rho}_u(x)=\ddot{u}$ if $\rho_u(x)\in \Phi_u$, and $\overline{\rho}_u(x)=\dot{u}$ if $\rho_u(x)\in B_u-\Phi_u$. Clearly $\overline{\rho}_u(u)=\ddot{u}$. Let $\Lambda_u=\{x\in V(M)|A_u(x)\neq\emptyset\}$. Further, for each distinct $x,y\in \Lambda_u$, let $\Lambda_u(x,y)=x\dot{u}y$ if either $x=u$, $y\in U_0$ or $y=u$, $x\in U_0$. Otherwise,
$$
\Lambda_u(x,y)=x\overrightarrow{\Upsilon}(x)\rho_u(x)\overline{\rho}_u(x)\Upsilon(u)\overline{\rho}_u(y)\rho_u(y)\overleftarrow{\Upsilon}(y)y.\qquad\qquad
$$

Let $\xi_1,...,\xi_f$ be the elements of $\Lambda_u$, occuring on $\overrightarrow{M}$ in a consecutive order with $\xi_1=u$. For each integer $i$ $(1\leq i\leq f)$, set
$$
M_i=\xi_i\overrightarrow{M}\xi_{i+1}, \quad \omega_i=|A_u(\xi_i)|+|A_u(\xi_{i+1})| \quad \mbox{(indices mod f)}.
$$

\textbf{Claim 1}. $\sum^f_{i=1}|M_i|\geq\sum^f_{i=1}\omega_i.$ \\

Proof. Since $M$ is extreme, we have $|M_i|\geq\omega_i$ for each $i\in \{2,...,f-1\}$. If $\Phi_u\cap V(\Upsilon(\xi_2))\neq\emptyset$ and $\Phi_u\cap V(\Upsilon(\xi_f))\neq\emptyset$, then clearly $|M_i|\geq\omega_i$ $(i=1,f)$ and we are done. Now let $\Phi_u\cap V(\Upsilon(\xi_2))=\emptyset$ and $\Phi_u\cap V(\Upsilon(\xi_f))=\emptyset$. Clearly $|M_i|\geq\omega_i-|A_u(u)|+1$ $(i=1,f)$. By the definition of $\Lambda_u$, $\dot{u}\xi_2\in E$ and $\dot{u}\xi_f\in E$. Since $u\in U_1$, we have $\Phi_u\cap V(\Upsilon(\xi_s))\neq\emptyset$ for some $3\leq s\leq f-1$. Then we can choose $i,j$ such that $2\leq i\leq s-1$ and $s\leq j\leq f-1$ with $|M_i|\geq\omega_i+|A_u(u)|-1$ and $|M_j|\geq\omega_j+|A_u(u)|-1$ and the result follows. Finally, because of the symmetry, we can suppose that $\Phi_u\cap V(\Upsilon(\xi_2))=\emptyset$ and $\Phi_u\cap V(\Upsilon(\xi_f))\neq\emptyset$. Clearly $|M_1|\geq\omega_1-|A_u(u)|+1$. By the definition of $\Lambda_u$, $\dot{u}\xi_2\in E$. Then we can choose $i$ such that $2\leq i\leq f-1$ with $|M_i|\geq\omega_i+|A_u(u)|-1$ and again the result follows. \quad $\Delta$\\

\textbf{Claim 2}. If $|\Upsilon(u)|\geq2$, then $\Phi_u\cap U_0=\emptyset$.\\

Proof. Suppose to the contrary and let $v\in \Phi_u\cap U_0$. Then replacing $\Upsilon(u)$ and $\Upsilon(v)$ by $u\overrightarrow{\Upsilon}(u)(\ddot{u})^-$ and $v\ddot{u}$, respectively, we can form a new $M_H$-spreading, contradicting the $(U_0)$-minimality of $\Upsilon$. \quad $\Delta$\\

By (2) and Claim 1,
$$
|M|=\sum^f_{i=1}|M_i|\geq\sum^f_{i=1}\omega_i=\sum^f_{i=1}(|A_u(\xi_i)|+|A_u(\xi_{i+1})|)
$$
$$
\qquad =2\sum^f_{i=1}|A_u(\xi_i)|=2\sum_{x\in V(M)}A_u(x)=2|\Phi_u\cup B_u|. \eqno(3)
$$

If $|\Upsilon(u)|\geq2$, then by Claim 2, $|\Phi_u\cup B_u|=\varphi_u+b_u$, which by (3) gives $|M|\geq2(\varphi_u+b_u)\geq\varphi_u+b_u+1$. Finally, if $|\Upsilon(u)|=1$, i.e. $\ddot{u}=\dot{u}$, then $|\Phi_u\cup B_u|=|B_u|+|\{u\}|=b_u+1=\varphi_u$ and again by (3), $|M|\geq2|\Phi_u\cup B_u|\geq2\varphi_u\geq\varphi_u+b_u+1$. \quad $\Delta$\\

\noindent \textbf{Proof of Lemma 3}. Put $L=u_1...u_m$. Let $\Lambda_u$, $\Lambda_u(x,y)$ and $\omega_i$ be as defined in proof of Lemma 2. Let $\xi_1,...,\xi_f$ be the elements of $\Lambda_u$ occuring on $\overrightarrow{L}$ in a consecutive order. Set
$$
\overrightarrow{M}^\prime=u_1\overrightarrow{L}\xi_1, \quad \overrightarrow{M}^{\prime\prime}=\xi_f\overrightarrow{L}u_m, \quad\overrightarrow{M}_i=\xi_i\overrightarrow{L}\xi_{i+1} \quad (i=1,...,f-1). 
$$

Let $G^\prime$ be the graph obtained from $G$ by adding an extra edge $u_mu_1$. Set $\overrightarrow{M}=u_1...u_mu_1$ and $M_f=\xi_f\overrightarrow{M}\xi_1$. Let $\Lambda^\prime_u(\xi_f,\xi_1)$ and $\Lambda^{\prime\prime}_u(\xi_f,\xi_1)$ be the paths obtained from $\Lambda_u(\xi_f,\xi_1)$ by deleting the first and the last edges, respectively. Since $L$ is extreme, $|M_i|\geq|\Lambda_u(\xi_i,\xi_{i+1})|$ $(i=1,...,f-1)$. As for $M_f$, observe that 
$$
|M^\prime|\geq |\Lambda^{\prime}_u(\xi_f,\xi_1)|=|\Lambda_u(\xi_f,\xi_1)|-1,
$$
$$
|M^{\prime\prime}|\geq |\Lambda^{\prime\prime}_u(\xi_f,\xi_1)|=|\Lambda_u(\xi_f,\xi_1)|-1,
$$
implying that
$$
|M_f|=|M^\prime|+|M^{\prime\prime}|+1\geq2|\Lambda_u(\xi_f,\xi_1)|-1\geq|\Lambda_u(\xi_f,\xi_1)|.
$$

So, $|M_i|\geq|\Lambda_u(\xi_i,\xi_{i+1})|$ for each $i\in\{1,...,f\}$. Further, for each $u\in U_1$, we can argue exactly as in proof of Lemma 2 to get $|L|=|M|-1\geq\varphi_u+b_u$.

Now let $u\in U^*$. By the definition, $\Phi_u\subseteq V(\Upsilon(u))$ and therefore, $|\Upsilon(u)|\geq|\Phi_u|=\varphi_u$. Since $L$ is extreme, $|L|\geq2(|B_u|+|\{u\}|)-2=2b_u$. Hence,
$$
|L|\geq |\Upsilon(u)|+\frac{1}{2}|\overrightarrow{L}|\geq\varphi_u+b_u. \qquad\qquad \Delta   \qquad\qquad\qquad\enspace
$$

\noindent\textbf{Proof of Theorem 1}. Let $M$ be a longest path in $G-C$ of length $\overline{p}$ with a $(U_0)$-minimal $M_C$-spreading $\Upsilon$. If $\overline{p}=-1$, i.e. $M$ is a Hamilton cycle, then $|C|\geq\delta+1=(\overline{p}+2)(\delta-\overline{p})$. Let $\overline{p}\geq0$. We claim that\\

\textbf{(a1)} \quad if $u\in U_0$ and $v\in \overline{U}_0$, then $\Phi_u\cap V(\Upsilon(v))\subseteq \{v,\dot{v}\}$,\\

\textbf{(a2)} \quad if $u\in U_0$, then $\varphi_u\leq \overline{p}+b^*_u$, \\
 
\textbf{(a3)} \quad if $v\in \overline{U}_0$, then $\varphi_u\leq \overline{p}-b_u$.\\

Let $u\in U_0$. If $v\in \overline{U}_0$, then to prove (a1) we can argue exactly as in proof of Claim 2 (see the proof of Lemma 2). The next claim follows immediately from (a1). To prove (a3), let $v\in \overline{U}_0$. Since $M$ is extreme, by Lemma 3, $\overline{p}\geq \varphi_u+b_u$ for each $u\in\overline{U}_0$, and (a3) follows. 

Observing that $\sum_{u\in U_0}b^*_u=\sum_{u\in \overline{U}_0}b_u$ and using (a2) and (a3), we get
$$
\sum_{u\in V(M)}\varphi_u\leq \overline{p}(\overline{p}+1)+\sum_{u\in U_0}b^*_u-\sum_{u\in \overline{U}_0}b_u=\overline{p}(\overline{p}+1).\quad\qquad
$$

Since $\Upsilon$ is extreme, we have $\psi_u=d(\ddot{u})-\varphi_u \geq\delta-\varphi_u$ for each $u\in V(M)$. By summing, we get
$$
\sum_{u\in V(M)}\psi_{u}=(\overline{p}+1)\delta-\sum_{u\in V(M)}\varphi_{u}\geq (\overline{p}+1)(\delta-\overline{p}).   \qquad \Delta\\
$$

\noindent\textbf{Proof of Theorem 2}. Let $M=u_1u_2...u_{\overline{c}}u_1$ be a longest cycle in $G-C$ of length $\overline{c}$ with an $(U_0,U^*)$-minimal $M_C$-spreading $\Upsilon$. Put
$$
U^*_1=\{u\in U^*|\varphi_u\leq \frac{1}{2}\overline{c}\}, \quad U^*_2=\{u\in U^*|\varphi_u\geq \frac{1}{2}(\overline{c}+1)\}.
$$

We claim that\\

\textbf{(b1)} \quad if $u\in U_0$ and $v\in \overline{U}_0$, then $\Phi_u\cap V(\Upsilon(v))\subseteq \{v,\dot{v}\}$,\\

\textbf{(b2)} \quad if $u\in U_0$, then $\varphi_u\leq \overline{c}-1+b^*_u$,\\

\textbf{(b3)} \quad if $u\in U_1$, then $\varphi_u\leq \overline{c}-1-b_u$,\\

\textbf{(b4)} \quad if $u\in U^*$, then $\varphi_u\leq \overline{c}-1-b_u+\varphi_u-\frac{1}{2}\overline{c}$,\\

\textbf{(b5)} \quad if $u\in U_1\cup U^*_1$, then $\varphi_u\leq \overline{c}-1-b_u$.\\

The proof of (b1) is very similar to proof of (a1) (see the proof of Theorem 1). The next claim follows immediately from (b1). By Lemma 2, $\overline{c}\geq\varphi_u+b_u+1$ for each $u\in U_1$ and (b3) follows. Since $M$ is extreme, $\overline{c}\geq 2(b_u+1)$ for each $u\in U^*$, which is equivalent to (b4).  Finally, (b5) follows from (b3) and (b4), immediately.

If $U^*_2=\emptyset$, then by (b2) and (b5), $\sum_u\varphi_u\leq \overline{c}(\overline{c}-1)$ and as in proof of Theorem 1, $|C|\geq (\overline{c}+1)(\delta-\overline{c}+1)$. Now let $U^*_2\neq\emptyset$. Choose $v\in U^*_2$ such that
$$
\varphi_v=\max_{u\in U^*_2}\{\varphi_u\}. 
$$

Then from (b4) we get\\

\textbf{(b6)} \quad if $u\in U^*_2$, then $\varphi_u\leq \overline{c}-1-b_u+\varphi_v-\frac{1}{2}\overline{c}$.\\

Using (b2),(b5),(b6) and observing that $\sum_{u\in U_0}b^*_u=\sum_{u\in \overline{U}_0}b_u$ and $|U_0|+|U_1\cup U^*_1|+|U^*_2|=\overline{c}$, we get
$$
\sum_u\varphi_{u}=\sum_{u\in U_0}\varphi_u+\sum_{u\in U_1\cup U^*_1}\varphi_u+\sum_{u\in U^*_2}\varphi_u\leq\overline{c}(\overline{c}-1)+|U^*_2|(\varphi_v-\frac{1}{2}\overline{c}). \eqno(4)
$$

By the definition, $\Phi_v\subseteq V(\Upsilon(v))$. Let $v_1,...,v_t$ be the elements of $\Phi^+_v$, occuring on $\overrightarrow{\Upsilon}(v)$ in a consecutive order with $v_t=\ddot{v}$. Clearly $t=|\Phi_v|=\varphi_v$. Put
$$
N(v_i)\cap V(\Upsilon)=\Phi^\prime_i, \quad  N(v_i)\cap V(C)=Z^\prime_i \quad (i=1,...,t).  
$$

If $\Phi^\prime_i\not\subseteq V(\Upsilon(v))$ for some $i\in\{1,...,t\}$, then replacing $\Upsilon(v)$ by $v\overrightarrow{\Upsilon}(v)v^-_i\ddot{v}\overleftarrow{\Upsilon}(v)v_i$, we form a new $M_C$-spreading, contradicting the minimality of $|U^*|$. So, we can assume that $\Phi^\prime_i\subseteq V(\Upsilon(v))$ $(i=1,...,t)$. Assume w.l.o.g. that 
$$
\max_i|\Phi^\prime_i|=|\Phi^\prime_t|=|\Phi_v|=\varphi_v=t.  \eqno(5)
$$

Since $\psi_{u_i}=d(u_i)-\varphi_{u_i}\geq\delta-\varphi_{u_i}$ $(i=1,...,\overline{c})$ and $|Z^\prime_i|=d(v_i)-|\Phi^\prime_i|\geq\delta-|\Phi^\prime_i|$ $(i=1,...,t-1)$, we have
$$
\sum^{\overline{c}}_{i=1}\Psi_{u_i}+\sum^{t-1}_{i=1}|Z^{\prime}_i|=\sum^{\overline{c}}_{i=1}(\delta-\varphi_{u_i})+\sum^{t-1}_{i=1}(\delta-|\Phi^{\prime}_{i}|)
$$
$$
=\delta(\overline{c}+t-1)-\sum^{\overline{c}}_{i=1}\varphi_{u_i}-\sum^{t-1}_{i=1}|\Phi^\prime_i|. \qquad\qquad\qquad\eqno(6)
$$ 

\textbf{Case 1}. $|U^*_2|=1$.\\

By (4), (5) and (6),
$$
\sum^{\overline{c}}_{i=1}\psi_{u_i}+\sum^{t-1}_{i=1}|Z^{\prime}_i|\geq (\overline{c}+t-1)\delta-\overline{c}(\overline{c}-1)-t+\frac{1}{2}\overline{c}-\sum^{t-1}_{i=1}t
$$
$$
=(\overline{c}+t-1)\delta-(\overline{c})^2-t^2+\frac{3}{2}\overline{c}.\qquad\qquad\qquad\qquad\qquad\qquad
$$

It follows, in particular, that
$$
\max_i\{\psi_{u_i},|Z^{\prime}_i|\}\geq\delta-\frac{(\overline{c})^2+t^2-\frac{3}{2}\overline{c}}{\overline{c}+t-1}\geq\delta-\frac{3}{2}\overline{c}+2.
$$

If $\delta\leq\overline{c}-1$, then clearly $|C|\geq(\overline{c}+1)(\delta-\overline{c}+1)$. Let $\delta\geq\overline{c}\geq t+1$. Applying Lemma 1 to $Q=\ddot{v}\overleftarrow{\Upsilon}(v)v\overrightarrow{M}v^-$, we get
$$
|C|\geq \sum^{\overline{c}}_{i=1}\psi_{u_i}+\sum^{t-1}_{i=1}|Z^{\prime}_i|+\max_i\{\psi_{u_i},|Z^{\prime}_i|\}\qquad\qquad\qquad\qquad\quad
$$
$$
\quad \qquad \geq(\overline{c}+1)(\delta-\overline{c}+1)+(t-1)(\delta-t-1)\geq(\overline{c}+1)(\delta-\overline{c}+1).
$$

\textbf{Case 2}. $|U^*_2|\geq2$.\\

Choose $w\in U^*_2-v$ such that $\varphi_v\geq\varphi_w\geq\varphi_u$ for each $u\in U^*_2-\{v,w\}$. Let $w_i, Z^{\prime\prime}_i, \Phi^{\prime\prime}_i$ $(i=1,...,r)$ be the analogs of $v_i$, $Z^{\prime}_i$ and $\Phi^{\prime}_i$ (defined for $\Upsilon(v)$) defined in this case for $\Upsilon(w)$ . As in (5), we can assume w.l.o.g. that $\max_i|\Phi^{\prime\prime}_i|=|\Phi^{\prime\prime}_r|=|\Phi_w|=\varphi_w=r$. Clearly $t+r=\varphi_v+\varphi_w\geq\overline{c}+1$. Then
$$
\sum^t_{i=1}|Z^{\prime}_i|+\sum^r_{i=1}|Z^{\prime\prime}_i|=\sum^t_{i=1}(d(v_i)-|\Phi^{\prime}_i|)+\sum^r_{i=1}(d(w_i)-|\Phi^{\prime\prime}_i|)
$$
$$
\geq\delta\varphi_v+\delta\varphi_w-\sum^t_{i=1}|\Phi^{\prime}_i|-\sum^r_{i=1}|\Phi^{\prime\prime}_{i}|\geq(t+r)\delta-t^2-r^2.\quad
$$

In particular, 
$$
\max_i\{|Z^{\prime}_i|,|Z^{\prime\prime}_i|\}\geq\delta-\frac{t^2+r^2}{t+r}.
$$

Applying Lemma 1 to $Q=\ddot{v}\overleftarrow{\Upsilon}(v)v\overrightarrow{M}w\overrightarrow{\Upsilon}(w)\ddot{w}$, we obtain
$$
|C|\geq\sum^t_{i=1}|Z^{\prime}_i|+\sum^r_{i=1}|Z^{\prime\prime}_i|+\max_i\{|Z^{\prime}_i|,|Z^{\prime\prime}_i|\}\qquad\qquad\qquad\qquad\qquad
$$
$$
\geq(t+r)\delta-t^2-r^2+\delta-\frac{t^2+r^2}{t+r}\qquad\qquad\qquad\qquad\quad\qquad
$$
$$
\quad\qquad\geq(\overline{c}+1)(\delta-\overline{c}+1)+\delta(t+r-\overline{c})+(\overline{c})^2-1-t^2-r^2-\frac{t^2+r^2}{t+r}.\thinspace\thinspace
$$

If $\delta\leq \overline{c}-1$, then clearly $|C|\geq(\overline{c}+1)(\delta-\overline{c}+1)$. Otherwise,
$$
|C|\geq(\overline{c}+1)(\delta-\overline{c}+1)+\overline{c}(t+r)-1-t^2-r^2-\frac{t^2+r^2}{t+r}
$$
$$
\thinspace \geq(\overline{c}+1)(\delta-\overline{c}+1)+(\overline{c}-1)(t+r)-t^2-r^2.\qquad
$$

Then we can obtain the desired result observing that 
$$
\overline{c}-1\geq\max\{t,r\}\geq\frac{t^2+r^2}{t+r}.  \qquad \Delta
$$

\end{document}